\documentclass[a4,11pt,leqno]{article}

\usepackage{amsmath} 
\usepackage{amssymb} 

\usepackage{latexsym}

\usepackage[OT2,T1]{fontenc}

\numberwithin{equation}{section}
\allowdisplaybreaks

\title{
On the Work of Benjamin Olinde\,Rodrigues\,(1795-1851)
\\
--- in particuler, on \lq\lq\,Expression of Spatial Motions\,\rq\rq ---}

\author{Takeshi HIRAI}

\date*{}

\begin{document}

\setcounter{page}{1}
\setcounter{section}{0}

\maketitle

{\small
{\bf Abstract.} 
This is a translation of Proceedings of 22\,th Symposium 
on History of Mathematics, Tsuda University 2011, on the work of Benjamin Olinde Rodrigues and on his life. His chef-d'oeuvre is the work on Euclidean motion group in 1840. He invented {\it Rodrigues expression}\, of rotation and gave explicit calculation formula for product of two rotations, which might be considered as a discovery of quaternion product rule ahead of Hamilton. We follow a new proof of \'E.\;Cartan in his book on {\it spineurs}\, in 1938 for Rodrigues formula, 
which was called as Euler-Olind-Rodrigues formula mistakenly. We add as Appendix important parts of Lecture Note on applications of quaternion. There are given 
description of rotational movements in Rodrigues expression and an interesting compact formula for time derivative of rotation, applicable in many purposes.
}
\footnote{Translation (with small changes) of 
Proceedings of 22\,th Symposium on History of Mathematics in 
 {\it Reports of Institute for Mathematics and Computer Sciences, Tsuda University}, 
{\bf 33}\,(2012), 59--79.}

{\normalsize 

\section{Motivation for research report and others}\; 

Why I came to pay attention to Benjamin Olinde Rodrigues (1795/10/6\,--\,1851/12/17) is because, in fact, I wanted to know the historical process of the quaternion of Hamilton. Actually, in the workshop "Non-Commutative Harmonic Analysis" at B\k{e}dlewo, Poland, August of 2007, I gave a talk on projective representations of complex reflection groups containing symmetric groups, and necessarily I discussed about the trilogy [Sch1, 1904]\,--\,[Sch3, 1911] of Schur's theory on spin representations of finite groups and noticed that, in the third paper [Sch3], there appeared substantially a triplet of matrices which is called later as {\it Pauli matrices}, found independently by Physicist Pauli in [Paul, 1927]. About this historical comment, Prof. Marek Bo\.zejko gave me a question \lq\lq{\it How about the case of quaternion of Hamilton ?}\;\rq\rq, which means that in the theory of quaternion in 19\,th century, \lq\lq\,Is there something like such triplet of matrices\;?\;''\; It might be possible, and mathematicians nowadays know generally that the unit ball $\boldsymbol{B}$ consisting of quaternions with norm 1 gives a double covering group of the 3-dimensional rotation group $SO(3)$. However, at that time, I lacked totally such kind of historical point of view, and couldn't answer his question. I felt deeply that I know nothing about Hamilton to answer this question and that I ought to study this sometime in future.  

\medskip

{\bf Note 1.} Three matrices appeared in p.198 of Schur's paper [Sch3, 1911] are
\begin{eqnarray}
\label{2012-01-27-1}
\;\;
F = 
\begin{pmatrix}
1&0\\ 0&1
\end{pmatrix}
,\quad 
A = 
\begin{pmatrix}
0&1\\ 1&0
\end{pmatrix}
,\quad  
B=
\begin{pmatrix}
0&1\\ 
-1&0
\end{pmatrix},\quad 
 C=
\begin{pmatrix}
1&0\\ 
0&-1
\end{pmatrix},
\\
\nonumber
\quad 
{\rm with}\qquad 
\left\{
\begin{array}{l}
A^2=F,\;\;B^2=-F,\;\;C^2=F,\;\; CBA=F,
\\
AB=-BA=-C,\; BC=-CB=-A,\;CA=-AC=B,\;
\end{array}\right.
\end{eqnarray}
and they are used, as fundamental ingredients, to write down explicitly doubly-valued projective representations of symmetric groups  $\mathfrak{S}_n$ and alternating group $\mathfrak{A}_n$.

On the other hand, three matrices appered in Pauli's paper [Paul, 1927] and later called as Pauli matrices, express rotation moment (spin) of electron, act on (${\boldsymbol C}^2$-valued) wave function $\psi$ as follows:\,\footnote{See also [Hir9].} 
\begin{eqnarray}
\label{2011-02-06-11}
&&
\boldsymbol{s}_x(\psi) = 
\begin{pmatrix}
0&1\\ 1&0
\end{pmatrix}\psi\,,\quad 
\boldsymbol{s}_y(\psi)=
\begin{pmatrix}
0&-i\\ i&0
\end{pmatrix}\psi\,,\quad
 \boldsymbol{s}_z(\psi)=
\begin{pmatrix}
1&0\\ 0&-1
\end{pmatrix}\psi\,,
\end{eqnarray}
with commutation relation\quad $
[s_x,s_y]=2is_z,\;[s_y,s_z]=2is_x,\;[s_z,s_x]=2is_y,\;i=\sqrt{-1}.
$
\vskip.2em
Further, \'Elie Cartan discussed in his paper [Car1, 1913], without introducing the terminology such as {\it spineur} ({\it spinor} in French), representations of $SU(2)$ double covering of rotation group $SO(3)$, actually containing spin representations together with linear representations, but did not appear there any triplet of matrices.

\vskip1em 

In the meantime, I was asked in 2008 from Prof.\;Satoshi Kawakami to give a Summer Intensive Course of one week for Mathematics Master Course of Nara Education University, and so I proposed as a subject \lq\lq{\it Quaternion, 3-dimensional Rotation Group and Introduction to Representation Theory of Groups}\rq\rq\, 
and begun to prepare a Lecture Note. In that occason, when I searched documents about quaternion in website, I found a paper [Alt2] {\it Hamilton, Rodrigues, and the Quaternion scandal}.\, What does it mean {\it scandal}, in such an academic situation ? I wondered and read it with a keen interest, then I noticed that the historical facts which I imagined myself until then, arround universal covering group of 3-dimensional rotation group, or parameter expressions of 3-dimensional rotations and so on, differ considerably from the true story. Thus I gradually went deeply into Rodrigues, the central figure of the present paper. As for Lecture Note above, I attach its important part with deep relation to this paper as an Appendix for more detailed commentary.

\section{Chronology for Benjamin Olind Rodrigues}

1789--1799: French Revolution,\quad 

1804--1814, 1815 (The One Hundred Days): Napoleon I，
 
1814--1815, 1815--1824: Burbon Restoration, Louis XVIII, 

1824--1830: Restoration, Charles X, 

1830--1848: Louis-Philippe I.
\vskip1em

\noindent
1795/10/06,  Benjamin Rodrigues was born. 
\\[1ex]
\indent
\quad 1807, Jews living in France were required to modify their family name. 

\quad 1808, Jews were required to add a name of a Christian Saint.  
\\[1ex]
\noindent
1811, Entrance examination for \'Ecole Polytechnique and \'Ecole Normale,

\hspace*{3ex} Rodrigues was ranked first (15 years old).
\vskip.7em
\noindent
Ranked first in the entrance examination, but he did not enter any of the above  \'Ecole. About this situation, I read several explanations: e.g., he was refused to enter because of Jews, or his oncre refused to pay his school expences, or next year he is also ranked first or second in the entrance exam, etc. Anyhow, for usual French people, at the same time of entering in one of these \'Ecoles, he is adopted by a government official, receive school expenses exemption, and a salary is paid, but it was not the case for Rodrigues. Finally he entered Universit\'e de Paris. 
\\[1.5ex]
\noindent
1815/06/23,  Soutien pour le doctorat (19 years old), 

\quad \`a la Facult\'e des Sciences de Universit\'e de Paris, sous la pr\'esidence de M.\;Lacrois, 

\quad Doyen de la Facult\'e.
\\
 1815/06/28, Awarded a doctorate in Mathematics 

\quad from Facult\'e des Sciences, Universit\'e de Paris.
 \vskip.7em
 \indent
\quad 1815\,$\sim$, After the 1815 Restoration the Catholic hierarchy took control of educational 

\quad and academic institutions, and Jewish people could not obtain any teaching position.
\\[1ex]
\noindent
 1816,\;\; The main part of the above th\`ese was published in 

\quad  {\it M\'emoire {\rm (*)} sur l'attraction des sph\'ero\"\i des, } 
\;{\sc PREMI\`ERE PARTIE.}   
{\it Formules

\quad g\'en\'erales pour l'attraction des corps quelquonques, et application de ces formules

\quad \`a la sph\`ere et aux ellipso\"\i des.}\; 
{\sc SECONDE PARTIE.} 
{\it Attraction des sph\'ero\"ides 

\quad infiniment peu diff\'erens d'une sph\`ere, et d\'evelopement g\'en\'erale de la fonction V.}

\quad ibid., {\bf 3}\;(1814--1816): pp.361--385, 1816.

\quad
 This paper contains {\it Rodrigues formula}\/ for Legendre polynomials\,:
$$P_n(x)=\dfrac{1}{2^n\,n!}\dfrac{d^n}{dx^n}(x^2-1)^n. $$
\vskip.2em

~ ----- between 1817--1837, no papers in Mathematics were published -----
\\[1ex]
1838, \;Three papers were published in Journal de Math\'ematiques pures et appliqu\'ees, 

\hspace*{3ex} vol.\;{\bf 3}\;(1938), pp.547--548,\; pp.549--549,\; pp.550--551. \hspace{8ex}(He was 42 years old)

\hspace*{3ex} The contents of 3 papers are the following:

\hspace{1ex}[R3-1] {\it The number of ways to decompose a convex polygon into triangles by diagonals} 

\hspace{1ex}[R3-2] {\it The number of ways to make product with $n$ factors}

\hspace{1ex}[R3-3] {\it Elementary and purely algebraic proof of development of binomial $(1+a)^x$ into 

\hspace*{7ex} power series with negative or frational powers}
\vskip.5em

Length of each paper is short, but contents are important. The first two gave compact 

and clear proofs for main results of long papers published very recently in journals.   
\\[1.5ex]
\noindent 
1839, \;A paper [R4-1] was pulished in {\it ibid.}, {\bf 4}(1839), 236--240. Its contents are 

\hspace{1ex}[R4-1] {\it Number of inversions in the order of products of permutations}
\vskip.5em

{\bf Contents.}\;  
The generating function for this number was given, and this result continued to have 
important influences until nowadays.\; 
For $\sigma\in\mathfrak{S}_n$, count the number of pairs $i<j$ such that  the inversion \;\lq\lq\,$\sigma(i)>\sigma(j)\,$'' occurs, and let  $N_n(k)$ be the number of $\sigma$ for which the number of inversins is just equal to $k$. Define the generating function of $N_n(k)$ as 
$$
R_n(q):={\sum}_kN_n(k)\,q^k.
$$
Rodriques gave a method of computing $R_n(q)$ inductively.

{\bf Note.}\; For this problem, later there were studies of Netto and MacMahon.  But, only\\ 
\hspace*{11ex} in 1970, Leonard Carlitz discovered this paper [R4-1]. 
\\[1.5ex]
\noindent
1840, \;A paper [R5-1] was pulished in {\it ibid.} vol.\;{\bf 5}, pp.380--440. Its subjects are 

 \hspace{1ex}[R5-1] {\it Eulidean Motions in the Space, in particular Rotations}
\vskip.5em

 This long paper gives very important results and can be said as his Chef d'Oeuvre.

The contents, containing substantial discovery of Quaternion, will be explaind in \S 3.  
\\[1.5ex]
\noindent
1843, \;Two papers were published in {\it ibid.}\,vol.\;{\bf 8}, 

\hspace*{3.5ex} [R8-1] pp.217--224,\; [R8-2] pp.225--234.  

\vskip1.2em

Cf. 1843/10/16, Didscovery of {\it Quaternion}\, by William Rowan Hamilton (1805--1865)\,:

\hspace*{16.5ex} by giving the fundamental formula \;\; $i^2=j^2=k^2=ijk=-1$.

\hspace*{3.3ex} 1843/10/17, Letter from Sir W.\,R.\;Hamilton to John T.\;Graves, Esq. 

\hspace*{5.2ex} {\it on Quaternions}\; (which is hand written)\quad [Later, in printed form, 1 line about 

\hspace*{5.2ex}  75 letters and total 167 lines, 4 pages,  containing footnotes, in Collected Works.]

\hspace*{3.3ex}  {\bf Note.} In these two days, he wrote, with a quill and ink, a letter of this volume  

\hspace*{5.2ex} and its fair copy for himself. This is a great concentration of a man of genius. 

\vskip1em
\noindent
1851/12/17,  Benjamin Olinde Rodrigues（56 years old）died.

\section{Chef-d'Oeuvre Paper\,[R5-1] in 1840} 

Paper [R5-1] by Olinde Rodrigues: {\it  Des lois g\'eom\'etriques qui r\'egissent les d\'eplacements d'un syst\`eme  
 solide dans l'espace, \,et la variation des coordonn\'ees provenant de ses 
 d\'eplacements consid\'er\'es ind\'ependamment des causes qui peuvent les produire.}

Translation of Title: 
{\it Geometric rules which govern the movements of a system of solid bodies in the space, and changes of coordinates coming from its displacements considered independently of the reason which produces them.} 
\vskip1.2em

{\bf The Style of Writing of this paper.} There is no independent Introducion. The paper is separated into parts with numbers from {\bf 1} to {\bf 33}, which we call here as {\it parts}\, (maybe {\it numeros} in French). Each {\it part}\, has only its number and no title. As general style, there are 18 Titles in italic under each of which a group of {\it parts}\, are gathered. But there exist several exceptions, for instance in some {\it parts}, there are one or two italic titles (something like as subsections in a section). Theorems are not separated from ground sentences as in modern style where theorems are numbered and their assertions are written in italic. In the middle of {\it part} 4 (or n$^o$\,4), there is a big italic title such as 
\begin{center}
{\it 
Th\'eor\`eme fondamental.
}
\end{center}
But this might be a title of a subnumero (like subsubsection in a subsection) and I found several of such italic titles. There are neither Propositions, Lemmas nor Diagrams.

\vskip1em

{\bf Contents. }\; He discussed very generally on displacements of solid bodies, that is, {\it Euclidean motion group}\, in modern language. Naturally there are two kind of motions, rotations and parallel translations. He treated the latter as a kind of {\it infinitesimally small rotations}\, with rotation axis situated in the perpendicular direction at the infinite long distances. This is the general idea (Id\'ee g\'en\'erale) throughout of this paper. Accordingly, he asserted (in n$^o$\,{\bf 1} and n$^o\,{\bf 2}$) that \lq\lq the properties of parallel translations are contained in the properties of rotations\,''. In his original expression, 

\begin{quotation}
Ainsi donc, toute {\it translation}\, d'un syst\`eme peut rigoureusement \^etre consid\'er\'e comme une rotation d'une amplitude infiniment petite autour d'un axe fixe infiniment \'eloign\'e et normal \`a la direction de cette translation.

On ne sera donc pas surpris de trouver ult\'erieurement toutes les propri\'et\'es des {\it translations}\, comprises dans celles des rotations, $\cdots\cdots$
\end{quotation}

Thus, as methods of discussions, repeatedly he used calculations using infinitely small pararell displacements  \,$\Delta x,\,\Delta y,\,\Delta z$, of the directions of $x$-axis and so on, and limit transitions from spheirical triangles to planar triangles.

\vskip1.2em

{\bf List of Titles and the corresponding numbers of \mbox{$\mathit parts$}\,:} 
\vskip.5em

{\it Id\'ee g\'en\'erale de la translation et de la rotation d'un syst\`eme solide.} 
\hspace*{3.7ex} n$^o$\,{\bf 1\;$\sim$\;2}

{\it 
Du d\'eplacement d'un syst\`eme d'un point fixe.} 
 \hspace*{26.3ex} n$^o$\,{\bf 3}
 
 \quad (Displacement fixing a point, or a rotation arround a fixing point.)
\vskip.5em

{\it 
Du d\'eplacement quelconque d'un syst\`eme solide dans l'espace.} 
\hspace*{9.7ex} n$^o$\,{\bf 4\;$\sim$\;7}

\vskip1ex
{\it 
De la composition des rotations successives d'un solide autour de

\hspace*{2.5ex}
deux axes convergents.
}
\hspace*{39ex} \hspace*{5.2ex} n$^o$\,{\bf 8}

\quad\; (Product of rotations arround two axis that intersect each other, or product of two 

\quad\, rotations arround the intersecting point.)

\noindent
This is one of highlight points of the paper where the composition of two rotations is treated. Here jumping over the usual product structure in the rotation group $SO(3)$, there appears the product structure in its universal covering (double covering) group Spin$(3)$. Thit is the {\it product formula in Rodrigues expression}\, of rotations. I quote his original sentences at this important critical point\,: 
\begin{quotation}
\noindent 
Telle est la diff\'erence caract\'eristique \`a signaler entre la composition des rotations et celle des translations successives. Il y a d'ailleurs entre des deux sortes de composition l'analogie qui existe entre les propri\'et\'es du triangle rectiligne et celles du triangle sph\'erique; et si l'on compare les translations parall\`eles aux trois c\^ot\'es d'un triangle rectiligne, aux sinus des demi-rotation accomplies autour des trois c\^ot\'es d'un angle tri\`edre, les valeurs des translations et celles de ces sinus seront \'egalement proportionnelles aux sinus des angles oppos\'es aux c\^ot\'es respectifs dans le triangle rectiligne et dans l'angle tri\`edre. 
\end{quotation}

{\bf (Explanation as I understand)}\; The first sentence, in response to the previous sentences, refers to \lq\lq difference between {\it rotation composition}\, and\, {\it translation composition}''. The following sentence explains how to calculate the compositions of two rotations arround the same center $O$. But it is difficult for me to translate this and relating portions of this n$^o$\,8 well into English correctly, so I will explain it in my own words.

On a unit sphere with center $O$, take two point $A$ and $B$ and put \,$\boldsymbol{n}_A=\overrightarrow{OA},\,\boldsymbol{n}_B=\overrightarrow{OB}$. Denote by $R(\phi_A \boldsymbol{n}_A)$ with \,$\phi_A\in {\boldsymbol R}$\, the rotation arround the unit vector $\boldsymbol{n}_A$ (as rotation axis) with the angle $\phi_A$ to the direction of right-handed screw. Then the assertion is 
\vskip1.2em 
{\bf Assertion 3.1.}\footnote{The title {\bf Assertion 3.1} is temporarily given here by me for convenience of quotation.} \;{\it The product $R(\phi_B \boldsymbol{n}_B)R(\phi_A \boldsymbol{n}_A)$ of two rotations is expressed as $R(\phi_C \boldsymbol{n}_C)$ with rotation axis $\boldsymbol{n}_C=\overrightarrow{OC}$ and rotation angle $\phi_C$ given as follows:\; 
rotate the plane $OAB$ arround $\boldsymbol{n}_A$ by angle $-\phi_A/2$, then we get a line (= big circle) on the sphere.\footnote{Note that if we rotate the plane $OAB$ by angle $\pi$ (a half of 2$\pi$), then it comes back to itself.} 
Similarly 
rotate the plane $OAB$ arround $\boldsymbol{n}_B$ by angle $\phi_B/2$, then we get another line on the sphere. Two lines intersect at a point $C$ (the nearest intersecting point), and we put inner angle at \,$C$ as $\pi-\phi_C/2$ (or outer angle $\phi_C/2$).
 }
\vskip1em

{\it Proof.}\; We draw several auxiliary lines (cf. Figure 3 in [Alt1]).   
Rotate the plane $OAB$ arround $\boldsymbol{n}_A$ by angle $\phi_A/2$, then we get a line on the other side of $\Delta ABC$, and rotate $OAB$ arround $\boldsymbol{n}_B$ by angle $-\phi_B/2$, we get another line, and they cross each other at a point $C'$.  $\Delta ABC'$ is a mirror image of $\Delta ABC$. Rotate $OAB$ arround $\boldsymbol{n}_A$ by angle $-\phi_A$, then we get a line on the other side of $\Delta ABC$, and rotate $OAB$ arround $\boldsymbol{n}_B$ by angle $\phi_B$, we get another line, and they cross each other at a point $A'$.  

{\bf 1)} \;$C$ is invariant under $R(\phi_B \boldsymbol{n}_B)R(\phi_A \boldsymbol{n}_A)$.

In fact, under the first rotation $R(\phi_A \boldsymbol{n}_A)$, $C$ is mapped to $C'$, and under the second $R(\phi_B \boldsymbol{n}_B)$, $C'$ is mapped back to $C$. 

{\bf 2)} \;$A$ is mapped to $A'$. \; 
In fact, under the first rotation, $A$ is mapped to $A$, and under the second, $A$ is mapped to $A'$. 

{\bf 3)} \;$\angle ACA'=\phi_C$. 

In fact, denote by $D$ the crossing point of $BC$ and $AA'$. Then $\Delta ACD$ and $\Delta DCA'$ are mutually mirror images of the other. So, $\angle ACA'=2\angle ACD=2(\phi_C/2)=\phi_C$. \hfill $\Box$

\vskip1.5ex

Thus Assertion 3.1 is proved and so $\sin(\phi_C/2)$ and $\cos(\phi_C/2)$ can be calculated using known formula in spherical trigonometry (cf. {\bf Note A3.1} in Appendix below). This gives a synthetic calculation method, and one can say that the calculation rule of {\it quaternion} has appeared substantially here  (cf.\,[Agn], [Alt1]). 
It seems that French mathematicians at the time had rather considerable background in spherical trigonometry.

\vskip1.2ex

{\it Composition des rotations infiniment petites.} 
 \hspace*{26.9ex} n$^o$\,{\bf 9}

\quad Here the translations are treated as infinitesimally small rotations arround an axis

\quad at infinity, and discuss their compositions (= products). 

\vskip1ex

{\it De la composition des rotations autour de deux axes parall\`eles.} 
 \hspace*{8.8ex} n$^o$\,{\bf 10\;$\sim$\;11}

\quad Composition of rotations arround two axis parallel to each other.

\vskip1ex
{\it De la composition des rotations autour d'axes fixes non 

\hspace*{2.5ex} 
convergents en nombre quelconques.} \hspace*{32.5ex} n$^o$\,{\bf 12}

\quad \;Composition of rotations arround non-intersecting two axis, very complicated.

\vskip1ex
{\it Examen du cas particulier des axes non convergents.}  
\hspace*{19ex} n$^o$\,{\bf 13} 

\quad \;Examination in the special case of two axis non-intersecting.

\vskip1ex

{\it De la composition des d\'eplacements successifs d'un syst\`eme 

\hspace*{2.5ex} 
 combin\'es de rotations et de translations.}
 \hspace*{27.5ex} n$^o$\,{\bf 14\;$\sim$\;15}

\quad On the composition of rotations and translations.

\vskip1ex

{\it \'Equation de l'axe central.}  
\hspace*{46.5ex} n$^o$\,{\bf 16} 

\quad Equations which determines the axis of the composed rotation. 

\vskip1ex

{\it Examen du cas des variations infiniment petites.} 
\hspace*{23.1ex} n$^o$\,{\bf 17\;$\sim$\;19} 

\quad Examination of the case of infinitely small displacements. 

\vskip1ex

{\it De la composition analytique des rotations autour d'axes non

\hspace*{2.5ex} 
convergents.}  
\hspace*{37ex}\hspace*{19.5ex} n$^o$\,{\bf 20} 

\quad Calculation formula for composition of rotations around two axis non-intersecting. 

\vskip.5em
{\it Composition des rotations successives autour de trois axes 

\hspace*{2.5ex} 
rectangulaires.}  
\hspace*{34.8ex}\hspace*{19.5ex} n$^o$\,{\bf 21}

\quad 
Composition of rotations arround three orthogonal axis (this corresponds so-called 

\quad 
Euler product of rotation).\footnote{By the way, Euler discussed the existence of axis for a rotation, but he didn't discuss 
Euler expression of a rotation, as I understand.}  

\vskip1ex

{\it De la composition des d\'eplacements infiniment petits successifs
 
\hspace*{2.5ex} 
d'un syst\`eme solide.} 
 \hspace*{29ex}\hspace*{19.5ex} n$^o$\,{\bf 22\;$\sim$\;23}

\vskip1ex
{\it Conditions d'\/{\rm \'equilibre} de plusieurs d\'eplacements successifs infiniment

\hspace*{2.5ex}
 petits.} 
 \hspace*{57ex} n$^o$\,{\bf 24}\;({\bf 25} missing)

\quad  Discussions on the realization of the state of {\it Equilibrium}, that is, the condition for 

\quad that, after successive displacements, the solid body comes back to the original 

\quad position. 

\vskip1ex

{\it Analogie de ces lois de composition et d'\'equilibre avec celles 
 de la composition et de 

\qquad\quad  l'\'equilibre des\, {\rm forces} appliqu\'ees \`a un syst\`eme invariable.}  
 \hspace*{8.5ex} n$^o$\,{\bf 26}

\quad  Discussions on the {\it striking analogy}\, (l'analogie frappante) between the above state 

\quad of {\it Equilibrium}\, and {\it Equilibrium state when force is applied}.

\vskip1ex

{\it De la d\'etermination des variations des coordon\'ees d'un syst\`eme 
 solid dues \`a un d\'e-

\qquad\quad placement quelconque de ce syst\`eme, 
analytiquement d\'eduites des conditions 

\qquad\quad de l'invariabilit\'e de ce syst\`eme.} 
\hspace*{33.3ex}  n$^o$\,{\bf 27\;$\sim$\;32}

\vskip1ex

{\sc Conclusion.}\quad --- {\it Loi g\'en\'erale de la Statistique.} 
\hspace*{22.4ex}  n$^o$\,{\bf 33}

\section{\'E.\,Cartan's Proof of Rodrigues formula for rotations}

\'Elie Cartan quoted in n$^o$\,{\bf 59.\;\,Repr\'esentation d'une rotation}, p.57, in his book [Car2], 1938, one of main results of Rodrigues\,\footnote{One of the main results of the paper [R5-1], 1840, of Benjamin Olinde Rodrigues. The author's name of this paper is written as Olinde Rodrigues.} as follows and gave a proof of his own. I quote the central part of the proof of Cartan:  
\begin{quotation}
La formule (3) permet de retrouver les formules d'Euler-Olinde-Rodrigues. Soit 
$L$ le vecteur unitaire\,\footnote{By definition, of the length 1.} port\'e sur l'axe de 
rotation et $\theta$\, l'angle de rotation\,; les deux vecteurs unitaires $A, B$ ont pour produit scalaire $\cos\tfrac{\theta}{2}$ et leur produit vectoriel $\tfrac{1}{2}(AB-BA)$ est \'egal \`a $iL\sin\tfrac{\theta}{2}$. On en d\'eduit
\begin{gather}
\label{2020-04-28-1}
\nonumber 
BA=\cos\tfrac{\theta}{2}-iL\sin\tfrac{\theta}{2},\quad 
AB=\cos\tfrac{\theta}{2}+iL\sin\tfrac{\theta}{2},
\end{gather}
\vskip-1em
\noindent
d'o\`u 
\vskip-2em
\begin{gather}
\nonumber
(5)\hspace{14ex}
X'=\Big(\cos\tfrac{\theta}{2}-iL\sin\tfrac{\theta}{2}\Big)X\Big(\cos\tfrac{\theta}{2}+iL\sin\tfrac{\theta}{2}\Big).
\hspace{15ex}
\end{gather}

Si l'on d\'esigne par $l_1,l_2,l_3$ les cosinus directeurs de $L$, les param\`etres d'Euler-Olinde-Rodrigues sont les quatre quantit\'es 
\begin{gather}
\label{}
\nonumber
\rho=\cos\tfrac{\theta}{2},\quad \lambda=l_1\sin\tfrac{\theta}{2},\quad \mu=l_2\sin\tfrac{\theta}{2},\quad \nu=l_3\sin\tfrac{\theta}{2},
\end{gather}
dont la somme des carr\'es est \'egal \`a 1.
\end{quotation} 

\noindent 
Here a matrix $X$ is associated to a vector $\overrightarrow{x}$ (as defined in n$^o$\,{\bf 55}) as 
\begin{gather}
\label{2020-04-29-1}
\nonumber
X=\begin{pmatrix} x_3 & x_1-ix_2 \\ x_1+ix_2 & -x_3 \end{pmatrix} 
\;\longleftrightarrow\;\overrightarrow{x}=\begin{pmatrix} x_1\\ x_2 \\ x_3 \end{pmatrix}\in{\boldsymbol R}^3, 
\end{gather}
and $\lambda E_2=\mbox{\rm diag}(\lambda,\lambda)$ is identified with a scalar $\lambda$, where $E_2$ is the identity matrix of order 2. Take a rotation $R(\theta\hspace{.2ex}\boldsymbol{l})$ arround a unit rotation axis $\boldsymbol{l}={}^t(l_1,l_2,l_3)$ and of rotation angle $\theta$. the matrix $L={\small \begin{pmatrix} l_3 & l_1-il_2 \\ l_1+il_2 & -l_3 \end{pmatrix}}$ is associated to the axis $\boldsymbol{l}$, and $X'$ is the image of $X$ under $R(\theta\hspace{.2ex}\boldsymbol{l})$. The matrices $A$ and $B$ are associated to unit vectors $\overrightarrow{a}$ and $\overrightarrow{b}$ respectively, chosen in such a way that \;$\tfrac{1}{2}(AB+BA)=\langle \overrightarrow{a},\overrightarrow{b}\rangle E_2=\cos\tfrac{\theta}{2}\,E_2,\;\tfrac{1}{2}(AB-BA)=iL\sin\tfrac{\theta}{2}$.  

Thus the formula (5) gives correctly Rodrigues expression of the rotation $R(\theta\hspace{.2ex}\boldsymbol{l})$ given in the paper [R5-1].\footnote{It is very much regrettable that Cartan did not give an exact reference to this paper, and he misunderstood the name of its author Olinde Rodrigues as names of two persons called Olinde and Rodrigues respectively (cf. {\bf Comment 5.1} below). }

\vskip1em

{\bf (Explanation of this quotation)}\; 
Let a matrix $X$ be associated to a vector $\overrightarrow{x}$ as above. 
Then we have\; $\det X =-(x_1^{\;2}+x_2^{\;2}+x_3^{\;2})E_2=-\|\overrightarrow{x}\|^2E_2$, where $\|\overrightarrow{x}\|$ denotes Euclidian norm of $\overrightarrow{x}$. Moreover\; $X^2=\|\overrightarrow{x}\|^2E_2$, and $A^2=E_2=1,\;A^{-1}=A$ \;for unit vector $\overrightarrow{a}$.  
Let $Y$ be associated to $\overrightarrow{y}$, then 
\begin{gather}
\label{2020-04-30-1}
\nonumber 
\tfrac{1}{2}(XY+YX)=\big\langle \overrightarrow{x},\overrightarrow{y}\big\rangle,\quad 
\tfrac{1}{2}(XY-YX)=i\, \overrightarrow{x}\!\wedge\! \overrightarrow{y},
\end{gather}
where \,$\overrightarrow{x}\!\wedge\!  \overrightarrow{y}$\, denotes the vector product of $\overrightarrow{x}$ and $\overrightarrow{y}$.\; Two vectors  $\overrightarrow{x}$ and $\overrightarrow{y}$ are perpendicular to each other if and only if \,$XY=-YX$, and in such a case \,$\overrightarrow{x}\!\wedge\!\overrightarrow{y}$ is called {\it bivecteur} (by Cartan) and is represented by $-i\,XY\,=-i\,\tfrac{1}{2}(XY-YX)$. 

In n$^o$\,{\bf 55}, a triplet $H_1,H_2,H_3$ of $2\times 2$ matrices is introduced as  
\begin{gather}
\label{2020-04-30-2}
\nonumber 
H_1=\begin{pmatrix} 0&1 \\ 1&0 \end{pmatrix},\quad 
H_2=\begin{pmatrix} 0&-i \\ i&0 \end{pmatrix},\quad
H_3=\begin{pmatrix} 1&0 \\ 0&-1 \end{pmatrix}, 
\end{gather}
and in n$^o$\,{\bf 57}, another triplet $I_j:=-iH_j\;(j=1,2,3)$ is introduced. Each of them has the following relations respectively 
\begin{gather}
\label{2020-04-30-3}
\nonumber 
H_j^{\;2} =1\;(j=1,2,3),\quad H_jH_k=-H_kH_j\;(j\ne k),\quad H_1H_2H_3=i,
\\
\nonumber 
I_j^{\;2} =-1\;(j=1,2,3),\;\;\quad I_jI_k=-I_kI_j\;(j\ne k),\;\;\quad I_1I_2I_3=-1. 
\end{gather} 
The matrix $X$ associated to $\overrightarrow{x}$\, is $X=x_1H_1+x_2H_2+x_3H_3$, and $\{I_1,I_2,I_3\}$ is a triplet satisfying the fundamental formula of {\it Quaternion} (this is remarked in  n$^o$\,{\bf 57, Relation avec la th\'eorie des quaternions}). Moreover $\{1, H_1,H_2,H_3,i,I_1,I_2,I_3\}$ gives a basis over ${\boldsymbol R}$ of $M(2,{\boldsymbol C})$ of full matrix algebra of order 2 over ${\boldsymbol C}$.

Now we come to explain the meaning of the above quotation. Consider a unit vector $\overrightarrow{a}$ and the reflection ({\bf sym\'etrie} in [Car2]) $\mathrm{Ref}(\overrightarrow{a})$ with respect to the hyperplain orthogonal to it, which is given as  
\begin{gather}
\label{2020-04-30-4}
\nonumber
\overrightarrow{x}'=\overrightarrow{x}-2\overrightarrow{a}\,\big\langle \overrightarrow{x},\overrightarrow{a}\big\rangle. 
\end{gather}
Translating this into the matrix form, we have\footnote{The fomula below is exacly \lq\,La formule (3)\,' at the top of the above quotation from [Car2], p.57.} 
\begin{gather}
\label{2020-04-30-5}
\mathrm{Ref}(\overrightarrow{a})X=X'=-AXA
\\
\nonumber
\mbox{\rm In fact,}\hspace{4ex}
X'=X-2A\,\frac{1}{2}(XA+AX)=X-AXA-A^2X=-AXA\qquad(\because\; A^2=1).
\hspace{12ex}
\end{gather}

Take another unit vector $\overrightarrow{b}$, then $\mathrm{Ref}(\overrightarrow{b})\mathrm{Ref}(\overrightarrow{a})X=BAXAB$. Also take a vector $\overrightarrow{y}$ perpendicular to $\overrightarrow{x}$ and take {\it bivecteur} $\overrightarrow{u}=\overrightarrow{x}\!\wedge\!\overrightarrow{y}$ which is represented by $U:=-i\,XY$. Then, under $\mathrm{Ref}(\overrightarrow{a})$, $XY$ is transformed to $X'Y'=A(XY)A=A(XY)A^{-1}$, and so\; 
$\mathrm{Ref}(\overrightarrow{b})\mathrm{Ref}(\overrightarrow{a})U=(BA)U(AB),\;AB=(BA)^{-1}$. 

Thus stated, we should come back to the fundamental principle of Cartan's idea for the proof of so-called\, {\it les param\`etres d'Euler-Olinde-Rodrigues}\, above.  In n$^o$\,{\bf 10.\;D\'ecomposition d'une rotation en un produit de sym\'etries} of [Car2], it is proved that, in the space of dimension $n$ over ${\boldsymbol R}$ or ${\boldsymbol C}$,  
\begin{quotation}
{\it Toute rotation est le produit d'un nomble pair $\leqslant n$ de sym\'etries.
}
\end{quotation}
Hence, in the case of $n=3$ over ${\boldsymbol R}$, every rotation is a product of two reflections ({\it sym\'etries}), that is, for any non-trivial rotation $R$, there exists two unit vectors $\overrightarrow{a}$ and $\overrightarrow{b}$, with the angle from $\overrightarrow{a}$ to $\overrightarrow{b}$ smaller than $\pi$, such that \,$R=\mathrm{Ref}(\overrightarrow{b})\mathrm{Ref}(\overrightarrow{a})$. 

The rotation axis $\overrightarrow{l}$ ($=\boldsymbol{l}$\, in our notation) of $R$ is a positive multiple of $\overrightarrow{a}\!\wedge\!\overrightarrow{b}$ and let the angle from $\overrightarrow{a}$ to $\overrightarrow{b}$ be $\tfrac{\theta}{2}$. Then $\big\langle \overrightarrow{a},\overrightarrow{b}\big\rangle =\cos\frac{\theta}{2}$\; or\; $\tfrac{1}{2}(AB+BA)= \cos\frac{\theta}{2}$\,, and\, $\tfrac{1}{2}(AB-BA)=i\,\overrightarrow{a}\!\wedge\!\overrightarrow{b}= \sin\frac{\theta}{2}\cdot iL$. 

If we check the movement of $\mathrm{Ref}(\overrightarrow{b})\mathrm{Ref}(\overrightarrow{a})$ on the 2-dimensional plane spanned by  $\overrightarrow{a}$ and $\overrightarrow{b}$, then it is exactly the movement of $R(\theta\hspace{.2ex}\boldsymbol{l})$, as we can see easily. 

Thus the assertion in the quotation above is newly proved by Cartan.

\section{Comments and Literatures}

\quad\; 
{\bf Comment 5.1.}\; For the middle name Olinde of Rodrigues, I have checked the lists of Christian Saints and the lists of traditional French boys' names, downlorded from website. Curiously enough, I couldn't find Olinde in these lists. I understood that this name is added by his father under the order of Christian Church arround 1808, but in some literature it is explained that this name Olinde, along with the second names of his brother and sisters, was taken by his father from literary works and the like. 
 Anyhow he signed to his mathematical papers as Olinde Rodrigues, not using Benjamin. The reason why, I cannot imagine, but this seems to work against him. For instance, \'Elie Cartan misunderstood it as two person's names, and in his book [Car2] used the terminology as {\it les formules} (and {\it les param\`etres}) {\it d'Euler-Olinde-Rodrigues}. 
My friend Prof.\,emeritus Michel Duflo helped me very much to search Rodrigues' papers, difficult to take copies. He wrote me that he didn't know there exists French name {\it Olinde} in that time.  

\vskip1.2em

 {\bf Comment 5.2.}\; About confusions on the first names of mathematicians, also there is the case of J.\;Schur (1875--1941) and I.\;Schur\;(Issai Schur). About 90 years later of the case of Benjamin Olind Rodrigues, it was the times when hidden (?) Jewish misanthropy turned into persecution by Nazi who held power. 
I cannot but shed tears in the latest years when Schur was becoming very unhappy because of the persecution (Cf. [Hir6]).

\vskip1.2em

{\bf Comment 5.3.}\; Some years ago, at about 2008, when I was checking systematically website files under the questioning title \lq mathematician Benjamin Rodreagues' or something similar, the considerable date and time were necessary while throwing away useless files to look for possibly valuable documents.
 In addition, 
there were many files which put wrong informations. 
For example, there appeared a file saying \lq the mathematician who wrote only one article during life' (it means [R5-1]), and after many files passing, \lq the mathematician who wrote only two articles during life', and so on. Also about his birth place and nationality, there are estimations \lq Spain or Portgal from the spelling of family name', and someone reported as \lq I found a family documents in a Portuguese 
ancient document house, so it's done' and so on. Finally I found a file reporting a third paper of Rodrigues, and I felt somethig really curious and so decided to study seriously the situation, in partricular, search how many papers are there of him etc. and wanted to collect copies of all of them. 

Once I asked Prof.\,M.\,Duflo to find out the above 3rd paper and so on in {\it Correspondence sur l'\'Ecole Imp\'eriale Polytechnique,\,{\bf 3}}. Then, together with copies of all papers of his in J.\;Math.\;Pures Appl., I could make a report on Rodrigues and sent him my draft. After that I again asked him to find another document [R4*] in {\it Bulletin Scientifique de la Soci\'et\'e Philomatique de Paris,} as shown below in a part-copy of his e-mail\,:   
\begin{quotation}
 Cher Takeshi,

Merci beaucoup pour ton int\'er\'essant et amusant texte sur Olinde. Je
suis content d'avoir pu t'aider \`a obtenir de la  documentation.
\begin{verbatim}
> Je te demandrai cette fois-ci encore de trouver l'article
> \bibitem[R04]{Rodr04} Olinde Rodrigues,
> Sur quelques propri\'et\'es des int\'egrales doubles et des rayons
> de courbure des surfaces,
> Bulletin Scientifique de la Soci\'et\'e Philomatique de Paris,
> pp.34-36, 1815. [Signed \lq P.' by Poisson.]
\end{verbatim}

Cela a \'et\'e difficile ! Un excellent exercice d'internet !

Je ne l'ai pas trouv\'e dans les biblioth\`eques parisiennes, y compris la
biblioth\`eque nationale, jussieu, ihp, ens... J'ai peut-\^etre mal cherch\'e,
je l'ai d\'ecouvert apr\`es, il s'appelle souvent Rodrigue (sans s) ce qui
rend les recherches difficiles.

Je ne l'ai pas trouv\'e sur le site web de la soci\'et\'e philomathique de paris.
En fait c'est amusant de consulter ce site web; je ne connaissais pas
l'\'existence de cette soci\'et\'e; $\cdots\cdots\cdots$

Puis j'ai pens\'e a books.google.com.  
C'est genial : j'ai mis "philomathique paris 1815" dans la boite de
recherches, et j'ai pu lire le livre --- et donc tu peux aussi le faire.
J'ignore dans quelle biblioth\`eque Google a scann\'e ce volume, probablement
une biblioth\`eque d'une universit\'e americaine. \;$\cdots\cdots$

\end{quotation}

\vskip.6em

{\bf Comment 5.4.} In this occasion (at the beginning of 2012), when I tryed again to look arround website putting \lq\,mathematician Benjamin Olinde Rodrigues\,' in the box of research, I was quite surprised that the situation has been completely changed from that of several years ago so that many files containing wrong data disappeared and the files, which should be appeared even then, such as [AlOr] etc. appear ranked high enough.

\vskip1.2em 

{\bf 5.5. Literatures.}\; 

Quoted from review [Dav] on the book [AlOr], AMS-LMS, 2005. 
\begin{quotation}
Rodrigues produced only 17 mathematical papers but wrote extensively about
social, economic, and political matters, on banking and on alleviating problems
of labor. From 1816 to 1837, Rodrigues produced no mathematical papers.
Between 1838 and 1845, he wrote eight, including one on transformation groups
that some consider his chef-d'oeuvre. Taken at face value, this is a remarkable achievement. How many of us could get back into mathematical
shape after doing something else (writing reviews or becoming a provost, say) for two decades\,? Perhaps Rodrigues was theorematizing all along
but didn't have the time to write up his findings properly. He left no personal papers, so we can't tell. We can safely conjecture, though, that he
kept abreast of the contents of the mathematical journals of the day. 
\end{quotation}

Quoted from \;[URL1]\; on Rodrigues.
\begin{quotation}
\noindent 
$\cdots\cdots$\;\; 
Rather in 1807 Jews living in France were required to modify their family names and in the following year they were required to add a name of French origin. At this point Olinde was added to Rodrigues name. \;$\cdots\cdots$
\end{quotation}

\vskip2em

{\Large\bf Appendix.\footnote{Lecture Note for Summer Concentrated Course for Mathematics Master Course (one week) of Nara Educational University, 2008}  
\;\; 

\quad Quaternion, 3-dimensional Rotation Group and

\quad Introduction to Theory of Reprentations of Groups}

\vskip.7em 

\setcounter{section}{2}

Section 1 is omitted. Section 2, Introduction and \S 2.1, are omitted. Below begins with \S 2.2. 
We add the character A to the top of section numbers and subsection numbers etc.

\vskip1.2em

{\Large\bf A2. Represent Complex Numbers and Quaternions by

\hspace*{4.2ex} means of $2\times 2$ Matrices}

\vskip.6em

\setcounter{subsection}{2}

{\large\bf A2.2. Complex $2\times 2$ matrices representing quaternions}

\vskip.3em

The total of {\it quaternion number} ${\boldsymbol q} = \alpha +\beta {\boldsymbol i} +\gamma{\boldsymbol j} +\delta {\boldsymbol k}\;\,(\alpha,\beta,\gamma,\delta \in {\boldsymbol R})$ gives a non-commutative number field, and we denote it by ${\boldsymbol H}$. Quaternion ${\boldsymbol q}$ with $\alpha=0$ is called {\it pure quaternion} and we denote by ${\boldsymbol H}_{\!-}$ their totality. Here ${\boldsymbol i},{\boldsymbol j},{\boldsymbol k}$ are imaginary units invented by W.R.\,Hamilton (1805--'65) satisfyng 
\begin{eqnarray}
\label{2008-02-23-1}
&&
\qquad 
{\boldsymbol i}^2 = {\boldsymbol j}^2 = {\boldsymbol k}^2 = \boldsymbol{ijk}=-1,\;\;\quad 
\\
&&
\label{2008-02-23-2}
\boldsymbol{ij}=-\boldsymbol{ji}=\boldsymbol{k},\quad \boldsymbol{jk}=-\boldsymbol{kj}=\boldsymbol{i},\quad \boldsymbol{ki}=-\boldsymbol{ik}=\boldsymbol{j},
\qquad \qquad 
\end{eqnarray}
and, together with 1, they form a basis of ${\boldsymbol H}$ over ${\boldsymbol R}$. Hamilton discovered on the way of morning walk on the 16th Oct.\;1843 that, with three imaginary units ${\boldsymbol i},{\boldsymbol j},{\boldsymbol k}$ satisfying the so-called {\it fundamental formula} (\ref{2008-02-23-1}), in the totality of $\boldsymbol{q} = \alpha +\beta \boldsymbol{i} +\gamma \boldsymbol{j} +\delta \boldsymbol{k}$, four arithmetic operations are possible. His paper [Ham1] appeared in 1843. Besides that, he wrote in the next day a detailed report about this dicovery to his friend J.T.\,Graves as a long letter, written with a quill and ink, which is reprinted in 1844. In the printed form of Collected Works of Hamilton, 4 pages, total 167 lines [Ham2].   

Quaternion ${\boldsymbol H}$ is a linear algebra over  ${\boldsymbol R}$ and can be immersed into $M(2,{\boldsymbol C})$. The immersion $\Psi$ is given as an linear extension of the correspondence, with $i=\sqrt{-1}\in {\boldsymbol C}$,   
\begin{eqnarray*}
\boldsymbol{i}\to I= \begin{pmatrix}
0 & -1 \\ 1 &0 
\end{pmatrix},
\quad 
\boldsymbol{j}\to J=
\begin{pmatrix}
0&i \\ i&0
\end{pmatrix},
\quad 
 \boldsymbol{k}\to K=
\begin{pmatrix}
-i&0 \\ 0&i
\end{pmatrix}
\end{eqnarray*}
so that for \,$\alpha,\alpha'\in {\boldsymbol R}$\, and \,$\boldsymbol{q,q}'\in{\boldsymbol H}$, 
\begin{eqnarray}
\label{2008-02-12-1}
\left\{
\begin{array}{rll}
\Psi(\alpha \boldsymbol{q}+\alpha'\boldsymbol{q}')&=&\alpha\Psi(\boldsymbol{q})+\alpha'\Psi(\boldsymbol{q}'), 
\\
\Psi(\boldsymbol{qq}') &=& \Psi(\boldsymbol{q})\Psi(\boldsymbol{q}').
\end{array}
\right.
\end{eqnarray}

The conjugate of $\boldsymbol{q} = \alpha+\beta \boldsymbol{i} +\gamma \boldsymbol{j} +\delta \boldsymbol{k}\in {\boldsymbol H}$ is defined as 
$
\overline{\boldsymbol{q}} = \alpha -\beta \boldsymbol{i} -\gamma \boldsymbol{j} -\delta \boldsymbol{k}
$
and the norm of $\boldsymbol{q}$ by 
$$
\|\boldsymbol{q}\|=\sqrt{\boldsymbol{q}\overline{\boldsymbol{q}}}=(\alpha^2+\beta^2+\gamma^2+\delta^2)^{1/2}.
$$
Then we have $\boldsymbol{q}^{-1} = \|\boldsymbol{q}\|^{-2}\overline{\boldsymbol{q}}$. 
\medskip  \medskip

{\bf Problem 2.2.1.} Prove that the system of relations (\ref{2008-02-23-2}) is equivalent with the system of relations (\ref{2008-02-23-1}). 
\vskip.8em

{\bf Problem 2.2.2.} Prove that the triplet of matrices  $\{I,J,K\}$ satisfies the similar relations as the triplet  $\{i,j,k\}$.

\vskip.8em

{\bf Problem 2.2.3.} Prove the formula (\ref{2008-02-12-1}). Also prove \;$
\det \Psi(\boldsymbol{q})= \|\boldsymbol{q}\|^2$,\, and \;$\overline{\boldsymbol{qq}'}= \overline{\boldsymbol{q}'}\;\overline{\boldsymbol{q}}$\; (the order of the product is inverted).

\bigskip

{\bf Note A2.2.4.} When Hamilton discovered quaternion on the way of morning walk, he was near to a bridge, and at that time he curved the so-called {\it fundamental formula} on a stone of the bridge. It is the formula (\ref{2008-02-23-1}). I read that even today some peoples of Department of Mathematics are used to take a morning walk to the bridge on the same date as the 16th October, the date of Grate Discovery.

\vskip1.2em
{\Large\bf A3.\quad Quaternion and rotation group $SO(3)$, Hamilton's \\
\hspace*{8.3ex}
 discrepancies}

\setcounter{section}{3}
\setcounter{equation}{0}

{\large\bf A3.1. \quad Expression of a 3-dimensional rotation}

\vskip.3em

After the discovery, Hamilton was pursuing applications of quaternion. One of the themes was the problem of describing a rotation in 3D Euclidean space $E^3$. 
As is explained in \S A2.1 (omitted), a rotation in 2D Euclidean plane $E^2\cong {\boldsymbol C}$ can be expressed in a simple way by multiplication of complex number with modulu 1, and so Hamilton was aiming for something similar to that with respect to quaternion. 

Let us take a bijective correspondence between the space ${\boldsymbol H}_{\!-}$ of pure quaternions and 3D Euclidean space  $E^3$ given as    
\begin{eqnarray}
\label{2008-02-20-21}
&&
{\boldsymbol H}_{\!-}\ni {\boldsymbol x} = x_1\boldsymbol{i}+x_2\boldsymbol{j}+x_3\boldsymbol{k}\;\longmapsto\;
x={}^t(x_1,x_2,x_3)=
\begin{pmatrix}
x_1\\x_2\\x_3
\end{pmatrix}
\in E^3.
\end{eqnarray}
Here $x$ is a vertical vector but we express it by a transposed horizontal vector to save space. The length of $x$ is given by 
$$
\|x\|=\sqrt{x_1^{\;2}+x_2^{\;2}+x_3^{\;2}}, 
$$
and so we put\, $\|\boldsymbol x\|:=\|x\|$. We can express the length preserving isomorphism (\ref{2008-02-20-21}) by a symbol as\, ${\boldsymbol H}_{\!-}\cong E^3$ (this modern symbol expression is powerful but didn't exist in the times of Hamilton !).

Recall that a rotation in $E^3$ fixing the origin is expressed by an orthogonal matrix $U=(u_{ij})_{i,j=1}^3\in SO(3)$ as 
 $x\to x'=Ux$, where, with $E_3$ the unit matrix of order 3, 
\begin{eqnarray}
\label{2008-02-20-31}
&&
SO(3)=\{U \in M(3,{\boldsymbol R});\;
U\,{}^tU= {}^tU\/U = E_3\,,\;\;\det U =1\}.
\end{eqnarray}

On the other hand, denote by ${\boldsymbol B}$ the unit ball of ${\boldsymbol H}$ given as \,${\boldsymbol B}:=\{\boldsymbol{a} \in {\boldsymbol H};\|\boldsymbol{a}\|=1\}$, then it is a group under multiplication. From the beginning, Hamilton estimated that, 
 as in the case of 2D rotation group $SO(2)$ and the torus group ${\boldsymbol T}^1:=\{z \in {\boldsymbol C};|z|=1\}$, 
\\[1.2ex]
\hspace*{10ex}{\it The group ${\boldsymbol B}$ should be isomorphic to the group $SO(3)$.} 
\vskip.9em 

He took as above the space of pure quaternions ${\boldsymbol H}_-$ as Euclidean space $E^3$ and look for ways of action of ${\boldsymbol B}$ on it. As we see, the simplest way of action is the left multiplication as \;$L(\boldsymbol{a}): {\boldsymbol H}_-\ni {\boldsymbol x}\,\mapsto\,{\boldsymbol a}{\boldsymbol x}\in {\boldsymbol H}_-\;\;(\boldsymbol{a}\in {\boldsymbol B})$. Alas ! for Hamilton. The image $L(\boldsymbol{a}){\boldsymbol x}={\boldsymbol a}{\boldsymbol x}$ belongs to ${\boldsymbol H}_-$ only when ${\boldsymbol a}$ is orthogonal to ${\boldsymbol x}$, that is, $\langle {\boldsymbol a},{\boldsymbol x}\rangle=\tfrac{1}{2}({\boldsymbol a}{\boldsymbol x}+{\boldsymbol x}{\boldsymbol a})=0$ (Cf. Problem 3.1 below). So the big and difficult problem for Hamilton was 
\vskip.8em

\noindent
{\it What kind of action of $\boldsymbol{a}\in {\boldsymbol B}$ on ${\boldsymbol H}_-\cong E^3$ gives an isomorphism from ${\boldsymbol B}$ onto $SO(3)$ ?}
\vskip1em

However, judging from the result, the above estimation of Hamilton was a misleading wrong estimate or a wrong button, from the beginning. Before 40 years old, He discovered quaterion， and after that he wrote several huge books e.g. [Ham3] and [Ham4], and tried to spread the theory of quaternion in the world, but it seems that the wrong buttons were stuck for his whole later life. 
In the essay [Alt2], Altmann wrote this situation in detail using some emotional terms such as {\it The sad truth}\, or {\it entirely unacceptable}\, or {\it Optical illusion}\, or {\it causing endless damage}\, etc. Still more expressed as 

\vskip.8em 
\quad {\it $\cdots\cdots$, and that Hamilton committed a serious error of judgement in basing
\\ 
\qquad his parametrization on the special case of the rectangular transformation.}\; 

\quad (this is the transformation appeared in Problem 3.1 below). 

\vskip.8em

{\bf Problem 3.1.}  Put ${\boldsymbol B}_-:={\boldsymbol B}\cap {\boldsymbol H}_{\!-}$ and express ${\boldsymbol a}\in {\boldsymbol B}$ as 
$$
{\boldsymbol a}=\cos \theta + \sin \theta \,{\boldsymbol w}\quad\mbox{\rm with}\;\;{\boldsymbol w} \in {\boldsymbol B}_-,\;\theta\in {\boldsymbol R}.
$$
Prove that, in case $\sin\theta\ne 0$, we have\; ${\boldsymbol a}{\boldsymbol x}\in {\boldsymbol H}_{\!-}$\; for an \;${\boldsymbol x}\in {\boldsymbol H}_{\!-}\cong E^3$\; if and only if \;$\langle {\boldsymbol x},{\boldsymbol w}\rangle=0$, that is, ${\boldsymbol x}\perp{\boldsymbol w}$ (perpendicular to each other).

Also prove that, in that case, the left multication $L({\boldsymbol a})=L(\cos \theta + \sin \theta \,{\boldsymbol w})$ induces on the hyperplane \;${\boldsymbol w}^\perp:=\{{\boldsymbol x}\in{\boldsymbol H}_{\!-}; {\boldsymbol x}\perp{\boldsymbol w}\}$\; a rotation of angle $\theta$ arround the origin.  

\vskip1.2em

Well now, {\it what is the correct expression of the rotation group by means of quaternion ?}\,  An answer has been given substantially in the paper [R5-1] of Rodrigues in 1840, but it is ignored historically until very recently, except an early comment by \'E.\,Cartan in [Car2].    

\vskip1.2em 
{\bf Lemma A3.2.}\; {\it For\; 
${\boldsymbol a}\in {\boldsymbol B}$, we have ${\boldsymbol a}^{-1}=\overline{{\boldsymbol a}},$ and the group ${\boldsymbol B}$ acts on ${\boldsymbol H}_{\!-}\cong E^3$ through 
\begin{eqnarray}
T({\boldsymbol a}):\quad {\boldsymbol H}_{\!-}\ni{\boldsymbol x}\,\longrightarrow\,{\boldsymbol x}' = {\boldsymbol a}{\boldsymbol x}{\boldsymbol a}^{-1}={\boldsymbol a}{\boldsymbol x}\overline{{\boldsymbol a}}\in{\boldsymbol H}_{\!-}, 
\end{eqnarray}
that is, \;$T({\boldsymbol a})T({\boldsymbol b})=T({\boldsymbol a}{\boldsymbol b})\;\;({\boldsymbol a},{\boldsymbol b}\in{\boldsymbol B})$.
}
\vskip1.2ex

{\it Proof.}\; Since\; $\overline{{\boldsymbol q}{\boldsymbol q}'}=\overline{{\boldsymbol q}'}\,\overline{{\boldsymbol q}}$, we have 
\;$\overline{{\boldsymbol x}'}=\overline{\overline{{\boldsymbol a}}}\,\overline{{\boldsymbol x}}\,\overline{{\boldsymbol a}}={\boldsymbol a}(-{\boldsymbol x})\overline{{\boldsymbol a}}=-{\boldsymbol x}'$,\; whence ${\boldsymbol x}'\in{\boldsymbol H}_{\!-}$.\; Moreover\; 
$T({\boldsymbol a})T({\boldsymbol b}){\boldsymbol x}={\boldsymbol a}({\boldsymbol b}{\boldsymbol x}\overline{{\boldsymbol b}})\overline{{\boldsymbol a}}
=({\boldsymbol a}{\boldsymbol b}){\boldsymbol x}(\overline{ {\boldsymbol a}{\boldsymbol b} })=T\big({\boldsymbol a}{\boldsymbol b} \big){\boldsymbol x}.$
\hfill
$\Box$
\vskip1.2em

{\bf Lemma A3.3.}\; 
{\it For a ${\boldsymbol w}\in {\boldsymbol B}_-={\boldsymbol H}_{\!-}\cap {\boldsymbol B},$ put \;$g_{\boldsymbol w}(\theta)=\cos\theta + \sin\theta\,{\boldsymbol w}\;(\theta\in \boldsymbol{R})$. Then $\theta\mapsto g_{\boldsymbol w}(\theta)$ is a one-parameter subgroup of ${\boldsymbol B}$, and \;$\dfrac{d\hspace{.2ex}g_{\boldsymbol w}(\theta)}{d\theta}\big|_{\theta=0}=\boldsymbol{w}$.
}

\vskip1em 

{\it Proof.}\; For $\theta, \theta'\in \boldsymbol{R}$, we have, from $\boldsymbol{w}^2=-1$,  
\begin{gather}
\label{2020-05-11-1}
\nonumber
g_{\boldsymbol w}(\theta)g_{\boldsymbol w}(\theta')=
\big(\cos(\theta)\cos(\theta')-\sin(\theta)\sin(\theta')\big)\hspace{10ex}
\\
\nonumber
\hspace{14ex}
+\big(\sin(\theta)\cos(\theta')+\cos(\theta)\sin(\theta')\big)\boldsymbol{w}
=g_{\boldsymbol w}(\theta+\theta').
\hspace{12ex}\Box
\end{gather}
\vskip.5ex

{\bf Problem 3.4.} \;For a ${\boldsymbol w}\in {\boldsymbol B}_-$, take ${\boldsymbol u},{\boldsymbol v}\in\boldsymbol{B}_-$ in such a way that \,${\boldsymbol u},{\boldsymbol v}\},\{{\boldsymbol w}$\, gives a right-handed orthonormal coordinate system. Then\; ${\boldsymbol u}{\boldsymbol v}={\boldsymbol w},\; 
{\boldsymbol v}{\boldsymbol w}={\boldsymbol u},\; 
{\boldsymbol w}{\boldsymbol u}={\boldsymbol v},$\; and for $g_{\boldsymbol w}(\theta)$, we have \ 
\begin{eqnarray}
\label{2020-05-11-12}
&&
\left\{
\begin{array}{l}
T(g_{\boldsymbol w}(\theta)){\boldsymbol u}=\quad\cos(2\theta){\boldsymbol u}+\sin(2\theta){\boldsymbol v},
\\
T(g_{\boldsymbol w}(\theta)){\boldsymbol v}=-\sin(2\theta){\boldsymbol u}+\cos(2\theta){\boldsymbol v},
\\
T(g_{\boldsymbol w}(\theta)){\boldsymbol w}=\quad{\boldsymbol w}.
\quad\;\,
\end{array}
\right.
\end{eqnarray}
The matrix expressin of $T(g_{\boldsymbol w}(\theta))$ with respect to the basis $\{{\boldsymbol u}, {\boldsymbol v},{\boldsymbol w}\}$ is 
\begin{eqnarray}
\label{2020-05-11-11}
\begin{pmatrix}
\cos(2\theta)& -\sin(2\theta)&0
\\
\sin(2\theta)&\;\;\cos(2\theta)&0
\\
0&0&1
\end{pmatrix}
.
\end{eqnarray}

\vskip1.2em

{\bf Theorem A3.5.}\; 
{\it The group ${\boldsymbol B}$ is a double covering and universal covering of rotation group $SO(3)$ and a covering map is given by $T({\boldsymbol a})\;({\boldsymbol a}\in{\boldsymbol B})$.
}
\vskip1em 

{\it Proof.}\; 
{\bf 1)}\; The map $T$ is surjective. In fact, it is known that any rotation $g\in SO(3)$ has a non-zero invariant vector $w\in E^3$ and so it is a rotation of of some angle $\phi$ arround $w$. Take the vector $\boldsymbol{w}\in {\boldsymbol H}_{\!-}\cong E^3$ corresponding to it and put $\theta=\phi/2$. Then, as seen from (\ref{2020-05-11-11}) in Problem 3.4, we have \;$T(\cos\theta+\sin\theta\,\boldsymbol{w})=g$.

{\bf 2)}\; The kernel of $T$ is $\{\pm 1\}\subset {\boldsymbol B}$. \;In fact, as seen from (\ref{2020-05-11-12}), $T(g_{\boldsymbol w}(\theta))=E_3$ (unit matrix) if and only if $2\theta\equiv 0\;(\mbox{\rm mod}\;2\pi)$, whence $\theta\equiv 0\;(\mbox{\rm mod}\;\pi)$ and so $g_{\boldsymbol w}(\theta)=\pm 1$. 

{\bf 3)}\; The unit ball ${\boldsymbol B}$ is topologically homeomorphic to 3-dimensional sphere $S^3$, and is simply connected. 
\hfill
$\Box$
\vskip1.2em

{\bf Problem 3.6.}\; Prove the following. 
 For $\boldsymbol{q} \in {\boldsymbol H}$, define $\exp \boldsymbol{q}$ by an abosolutely convergent infinite series as 
\begin{eqnarray}
\exp \boldsymbol{q} = \sum_{n=0}^\infty \frac{\boldsymbol{q}^n}{n!}=1+\boldsymbol{q}+\frac{\boldsymbol{q}^2}{2!}+\cdots.
\end{eqnarray}
Then, for 
$ {\boldsymbol w}\in {\boldsymbol B}_-$,\;  
$
\exp( \theta {\boldsymbol w}) = \cos \theta +  \sin\theta\,{\boldsymbol w}=g_{\boldsymbol w}(\theta)\;\; (\theta \in {\boldsymbol R}).
$

\vskip1.2em

The most important point of above dicussions is that, under the correspondence\;\begin{gather}
\label{2020-05-12-1}
{\boldsymbol H}_-\ni\theta\boldsymbol{w}\;\to\;\exp(\theta {\boldsymbol w})\in{\boldsymbol B}\;\to\; 
T\big(\exp(\theta {\boldsymbol w})\big)\in SO(3) 
\end{gather}
with \,$\theta\in\boldsymbol{R},{\boldsymbol w}\in {\boldsymbol B}_-$, the angle of rotation is doubled as $\theta\to2\theta$ as is shown in (\ref{2020-05-11-11}).
This means that $\boldsymbol{a}=\exp(\theta {\boldsymbol w})$ and $-\boldsymbol{a}=\exp\big((\theta+\pi){\boldsymbol w}\big)$ have the same image \;$T({\boldsymbol a})=T(-{\boldsymbol a})$, and the map \,$T: {\boldsymbol B}\to SO(3)$\, is a 2:1 correspondence. Hamilton seems to have been insisting particularly, with the great pioneer's stubbornness, to obtain 1:1 correspondence.
According to some biographies, Hamilton became eventually to suffer from excessive alcohol intake [Bell].

 Altmann points out, as one of the reasons, a serious psychological distress in this rotation expression problem. Looking at the cause of his worries, from the present age of mathematical standard, I can suspect that, in this problem there are two different objects such as   

1)\quad object which operate on something\quad(operators), 

2)\quad object to be operated\quad(operands), 
\\
however they both are the same quaternion and might be confused mutually or might not be clearly distinguished.  

Dear readers ! You may not 
feel much sympathy to Hamilton's serious anxiety, when reading this explanation. However it is because firstly you have been taught already under a modern mathematical basic training, and secondly here the author (Hirai) have chosen adequate notation in such a way that 

for the object 1), the characters such as ${\boldsymbol a},{\boldsymbol b}, g_{\boldsymbol w}$ etc., 

for the object 2), the characters such as ${\boldsymbol x},{\boldsymbol y}$ etc.\\
Thus, because of the hints drawn carefully for the reader, your understanding is unconsciously guided in the right direction.
\vskip1.2em

I have already mentioned the misunderstandings that Hamilton had. For his life and also about quaternion, somewhat ironic story telling can be found in website  

 http://members.fortunecity.com/jonhays/clifhistory.htm 
 \\
and the author of this website MR.\,jonhays noted that at age 17 he read about  Hamilton in {\it Men of Mathematics}\, by Eric Temple Bell\;(1883--1960).\footnote{Chapter 19, {\it An Irish Tragedy}, pp.340--361.} 

\vskip1em

{\bf Problem 3.7 (formula for calculation).} For $\boldsymbol{u, v} \in \boldsymbol{H}_{\!-}$, the corresponding elements in $E^3$ are denoted by $u= {}^t(u_1,u_2,u_3), v={}^t(v_1,v_2,v_3)$, and their inner product is defined as  
\;$\langle u,v\rangle := u_1v_1+u_2v_2+u_3v_3$ and denoted by  
$ \boldsymbol{u}\cdot\boldsymbol{v} =\langle u,v\rangle$. Then, prove the following formula:  
\begin{eqnarray}
\label{2008-02-25-1}
&&
\qquad
\boldsymbol{uv}= -\boldsymbol{u}\cdot\boldsymbol{v} + \boldsymbol{u} \times \boldsymbol{v},
\\
\nonumber
\mbox{\rm where}\qquad 
&&
\boldsymbol{u} \times \boldsymbol{v}\in \boldsymbol{H}_{\!-}, \qquad 
\boldsymbol{u} \times \boldsymbol{v} \;:=\; 
\begin{vmatrix}
\boldsymbol{i} & u_1&v_1 \\
\boldsymbol{j} & u_2 & v_2 \\
\boldsymbol{k} & u_3 &v_3
\end{vmatrix}\,.
\hspace{21ex} 
\end{eqnarray}

{\bf Problem 3.8.}  \;Prove the following formula: 
\begin{eqnarray}
\label{2008-03-04-1} 
\boldsymbol{u} \times \boldsymbol{v} = - \boldsymbol{v} \times \boldsymbol{u}, \quad 
(\boldsymbol{u} \times \boldsymbol{v}) \perp \boldsymbol{u}, \quad (\boldsymbol{u} \times \boldsymbol{v}) \perp \boldsymbol{v}.
\qquad 
\end{eqnarray}

\vskip.5em

{\large\bf A3.2.\; Rodrigues expression of rotation and product formula}

\medskip 

{\bf A3.2.1. \;Rodrighues parameter 
 $\theta \boldsymbol{w}\in \boldsymbol{H}_{\!-}\; 
(\theta\in {\boldsymbol R},{\boldsymbol w}\in {\boldsymbol B}_{-})$} 

 It was proved by Euler that any rotation $\rho$\, in 3D Euclidean space $E^3$ arround the origin has necessarily a rotation axis. Let $w\in E^3, \|w\|=1,$ be the axis and $\theta$ the angle of $\rho$ arround the axis $w$ in right-handed screw rotation. Then, as is shown above, $\rho$ is expressed in (\ref{2020-05-12-1}) as 
\begin{eqnarray}
\label{2008-02-21-1}
\rho=R(\theta\boldsymbol{w}):=T\big((\exp(\tfrac{1}{2}\theta{\boldsymbol w})\big)\in SO(3),
\end{eqnarray}
where $\boldsymbol{w}\in {\boldsymbol H}_{\!-}\cong E^3$ corresponds to $w$. 
 We call this expression as {\bf Rodrigues expression}, and \;$(\theta, \boldsymbol{w})\in {\boldsymbol R}\times {\boldsymbol B}_-$\; or \;$\theta\boldsymbol{w}\in \boldsymbol{H}_{\!-}$\; as {\bf Rodrigues parameter} of rotation $\rho$. 

This expression is, unlike the expression by means of Euler angles (cf. \S A3.3 below), the parameters are seamless, and locally univalent but globally multivalent. If one wishes to make completely univalent and put some restriction on $\theta {\boldsymbol w}$, there appears inevitably some breaks in the parameter. So that, as parameter space, it is natural to take the whole space ${\boldsymbol H}_{\!-}$ and enjoy the advantage of capability of describing smoothly multi-rotations of machines or airplanes etc.  

\vskip1em

{\bf A3.2.2. \;Product formula for two rotations}

The main contribution of Rodrigues is the description of the product of two rotations \;$R(\theta {\boldsymbol w})R(\theta' {\boldsymbol w}')=R(\theta'' {\boldsymbol w}'')$. \,To describe Rodrigues parameters $\theta'' {\boldsymbol w}''$ from \,$\theta {\boldsymbol w}$ and $\theta' {\boldsymbol w}'$, we calculate it according to the quaternion product rule in ${\boldsymbol H}$ as  
\begin{equation}
\label{2008-02-20-33}
\nonumber
{\textstyle 
\big(\cos({\frac{1}{2}}\theta)+\sin(\frac{1}{2}\theta){\boldsymbol w}\big)
\big(\cos(\frac{1}{2}\theta')+\sin(\frac{1}{2}\theta'){\boldsymbol w}'\big)
=
\cos(\frac{1}{2}\theta'')+\sin(\frac{1}{2}\theta''){\boldsymbol w}''.
}
\end{equation}
It gives us the so-called {\bf Rodrigues formula} in [R5-1] in our notations as
\begin{eqnarray}
\label{2008-02-29-1}
\hspace*{10ex}
\left\{
\begin{array}{rll}
\cos(\tfrac{1}{2}\theta'')
&\!\!=\!\!&\!
\cos(\tfrac{1}{2}\theta)\cos(\tfrac{1}{2}\theta') 
- \sin(\tfrac{1}{2}\theta)\sin(\tfrac{1}{2}\theta'){\boldsymbol w}\cdot {\boldsymbol w}',
\\[1ex]
\sin(\tfrac{1}{2}\theta''){\boldsymbol w}''
&\!\!=\!\!&\!
 \cos(\tfrac{1}{2}\theta)\sin(\tfrac{1}{2}\theta'){\boldsymbol w}'
\!+\!\sin(\tfrac{1}{2}\theta)\cos(\tfrac{1}{2}\theta'){\boldsymbol w}+
\\[1ex]
&& 
\hspace*{18.3ex}
+\sin(\tfrac{1}{2}\theta)\sin(\tfrac{1}{2}\theta'){\boldsymbol w}\!\times\!{\boldsymbol w}'.
\end{array}
\right.
\end{eqnarray}
Originally he induced his product formula (equivalent to the above) from some formulas for spherical trigonometric functions, 
and thus we may say that Rodrigues substantially gave the product rule for quaternion, in advance of Hamilton. 
\vskip1em

{\bf Note A3.9.} 
\;When $\theta,\theta'$ are both small, we can obtain the first approximation from the second equation above as 
\;$\theta'' {\boldsymbol w}''\doteqdot\theta'\,{\boldsymbol w}'+\theta\,{\boldsymbol w}.$\; Furthermore if 
$\boldsymbol{w},\boldsymbol{w}'$ are near to each other, then the first approximation is 
\begin{eqnarray}
\label{2008-02-29-11}
\theta'' \doteqdot \theta+\theta',\;\;\qquad  
\boldsymbol{w}'' \doteqdot\tfrac{1}{2}(\boldsymbol{w}+\boldsymbol{w}').
\qquad 
\end{eqnarray}
Note that, in the case of Euler angle expression, there does not exist such an approximation.  

\vskip1em

{\bf Note A3.10.} \;Basic formulas for spherical trigonometry are given as follows. 
For a spherical triangle $ABC$, let the interior angles be $\alpha, \beta, \gamma$,\, the opposite side lengths be $a, b, c$, the area be $S$, and the radius of the sphere be $\rho$. Then a formula {\it Spherical Excess} is 
\begin{eqnarray}
\label{2008-08-23-1}
&&
\alpha+\beta+\gamma -\pi=S/\rho^2 > 0,
\end{eqnarray}
\nonumber
Below, put radius \,$\rho=1$, then 
\begin{eqnarray}
\label{2008-08-23-2}
\mbox{\rm $\sin$ formula}
&&
\quad \frac{\sin a}{\sin \alpha} = \frac{\sin b}{\sin \beta} =\frac{\sin c}{\sin \gamma}, 
\\
\label{2008-08-23-3}
\mbox{\rm $\cos$ formula}
&&
\left\{
\begin{array}{l}
\cos a = \cos b\,\cos c + \sin b\,\cos c\,\cos\alpha, 
\\ 
\cos b = \cos c\,\cos a + \sin c\,\cos a\,\cos\beta, 
\\
\cos c = \cos a\,\cos b + \sin a\,\cos b\,\cos\gamma, 
\end{array}
\right.
\\
\label{2008-08-23-4}
\mbox{\rm $\cos$ formula}
&&
\left\{
\begin{array}{l}
\cos \alpha = -\cos \beta\, \cos \gamma + \sin\beta\,\sin \gamma\,\cos a, 
\\
\cos \beta = -\cos \gamma\, \cos \alpha + \sin\gamma\,\sin \alpha\,\cos b, 
\\
\cos \gamma = -\cos \alpha\, \cos \beta + \sin\alpha\,\sin \beta\,\cos c, 
\end{array}
\right.
\\
\label{2008-08-23-5}
\hspace*{5ex}\mbox{\rm $\sin\cos$ formula}
&&
\left\{
\begin{array}{l}
\sin a\,\cos\beta = \cos b\,\sin c - \sin b\, \cos c\, \cos \alpha, 
\\
\sin b\,\cos\gamma = \cos c\,\sin a - \sin c\, \cos a\, \cos \beta, 
\\
\sin c\,\cos\alpha = \cos a\,\sin b - \sin a\, \cos b\, \cos \gamma. 
\end{array}
\right.
\end{eqnarray}
In (\ref{2008-08-23-2})--(\ref{2008-08-23-5}), the radius $\rho=1$, and so the side length is equal to the angle measured in radians. For example, the length $c$ of $\overline{AB}$ is equal to \;$\angle AOB$ ($O$ is the center of sphere). The first formula in (\ref{2008-02-29-1}) comes from the 3rd formula in (\ref{2008-08-23-4}). In fact, 
$$
\alpha = \tfrac{1}{2}\theta,\;\;\beta=\tfrac{1}{2}\theta',\;\;
\gamma = \pi - \tfrac{1}{2}\theta'', \;\;c= \angle AOB,\;\;\cos c = \boldsymbol{w} \cdot \boldsymbol{w}'. 
$$

\vskip.7em

{\large\bf A3.3.\; Expression of 3D rotation by means of Euler angles} 

Let  $\{\boldsymbol{e}_1,\boldsymbol{e}_2,\boldsymbol{e}_3\}$ be an orthonormal system giving coordinate in $E^3$ as $x_1\boldsymbol{e}_1+x_2\boldsymbol{e}_2+x_3\boldsymbol{e}_3\;\leftrightarrow\; x={}^t(x_1,x_2,x_3)$. Denote by $g_1(\theta)$ a rotation arround $\boldsymbol{e}_1$ at right-handed screw angle $\theta$, and similarly $g_2(\theta),g_3(\theta)$ for $\boldsymbol{e}_2,\boldsymbol{e}_3$ respectively, then 
\begin{eqnarray}
\label{2008-02-27-51}
\nonumber
\;\;
g_1(\theta)\!=\! 
\begin{pmatrix}
1\!&\!0\!&\!0\\
0\!&\!\cos \theta\!&\!-\sin\theta \\
0\!&\!\sin\theta\!&\!\cos\theta\end{pmatrix}
,\,
g_2(\theta)\!=\!
\begin{pmatrix}
\cos \theta\!&\!0\!&\!\sin\theta \\
0\!&\!1\!&\!0\\
-\sin\theta\!&\!0\!&\!\cos\theta\end{pmatrix}
,\,
g_3(\theta)\!=\!  
\begin{pmatrix}
\cos \theta\!&\!-\sin\theta\!&\!0 \\
\sin\theta\!&\!\cos\theta\!&\!0\\
0\!&\!0\!&\!1
\end{pmatrix}
.
\qquad
\end{eqnarray}
A rotation $\rho$ of $E^3$ or $\rho\in SO(3)$ can be expressed as 
\begin{eqnarray}
\label{2008-02-27-55}
\nonumber
\;\;\rho
\!&=&\! 
g_3(\varphi)g_2(\theta)g_3(\psi) \hspace{20ex}(-\pi <\varphi,\psi \le \pi,\; 0\le \theta \le \pi)
\\
\nonumber
\!\!\!\!\!\!\!\!\!\!\!\!\!
&=&\!\!\!   
\begin{pmatrix}
\cos \varphi\!&\!-\sin\varphi\!&\!0 \\
\sin\varphi\!&\!\cos\varphi\!&\!0\\
0\!&\!0\!&\!1
\end{pmatrix} 
\begin{pmatrix}
\cos \theta\!&\!0\!&\!\sin\theta \\
0\!&\!1\!&\!0\\
-\sin\theta\!&\!0\!&\!\cos\theta\end{pmatrix}
\begin{pmatrix}
\cos \psi\!&\!-\sin\psi\!&\!0 \\
\sin\psi\!&\!\cos\psi\!&\!0\\
0\!&\!0\!&\!1
\end{pmatrix} 
\\
\label{2008-02-27-53}
\nonumber
\!\!\!\!\!\!\!\!\!\!\!\!\!
&=&\!\!\!
\begin{pmatrix}
\cos \varphi\cos\theta\cos\psi-\sin\varphi\sin\psi 
&-\sin\varphi\cos\theta\cos\psi-\cos\varphi\sin\psi&\sin\theta\cos\psi \\
\cos\varphi\cos\theta\sin\psi+\sin\varphi\cos\psi&-\sin\varphi\cos\theta\sin\psi+\cos\varphi\cos\psi&\sin\theta\sin\psi\\
-\cos\varphi\sin\theta&\sin\varphi\sin\theta&\cos\theta
\end{pmatrix} 
\!.
\qquad\qquad
\end{eqnarray}

\vskip.3em

When we apply Euler angle expression to calculations such as product of two rotations, we encounter immediately some dificulties, for instances: 

\medskip

{\bf 1)} The computational load is heavy. In fact, to calculate the product of two rotations, first calculate the product of 6 matrices of Euler angles, then decompose the product into 3 Euler angle components again. 
\medskip

{\bf 2)} It is not possible to evaluate calculation errors. Even if the change of rotations are small enough, the deviations between their Euler angles can be very big, sometimes there would be jumps. 
\medskip

For more details, you can read Altmann's text book [Alt1].

\vskip2em

{\Large\bf A4. Applications of Rodrigues expression, Time derivative}

\setcounter{section}{4}
\setcounter{equation}{0}
\vskip.5em

Recently, in many directions, Rodrigues expression of rotation is applied adequately.  Its parameter is given by $\theta\boldsymbol{w}\in \boldsymbol{H}_{\!-}\;(\theta\in \boldsymbol{R},\,\boldsymbol{w}\in \boldsymbol{B}_-)$ and the rotation $R(\theta\boldsymbol{w})=T\big((\exp(\tfrac{1}{2}\theta\boldsymbol{w})\big)$ acts on $\boldsymbol{x}\in \boldsymbol{H}_{\!-}\cong E^3$ as 
\begin{eqnarray}
\label{2008-02-28-1} 
&&
\quad 
\boldsymbol{H}_{\!-}\ni \boldsymbol{x} \;\to\; \exp(\tfrac{1}{2}\theta\boldsymbol{w})\,\boldsymbol{x}\,\exp(\tfrac{1}{2}\theta\boldsymbol{w})^{-1}
 \\
 \nonumber
 &&
 \hspace*{7ex}
= \big(\cos(\tfrac{1}{2}\theta)+\sin(\tfrac{1}{2}\theta)\boldsymbol{w}\big)\,
\boldsymbol{x}\,\big(\cos(\tfrac{1}{2}\theta)-\sin(\tfrac{1}{2}\theta)\boldsymbol{w}\big)\in\boldsymbol{H}_{\!-}\cong E^3,
\qquad 
\end{eqnarray}

\vskip.3em

{\bf A4.1. Examples of application in various fields}

In geophysics, it is important to describe rotations for problems such as\,:
\medskip

$\bullet$\; In plate tectonics, describe plates movement according to geological time, 
 by means of a rotation which leaves the earth center invariant.
\vskip.4em

$\bullet$\; In geodesy, describe the relationship of the inertial coordinate system of the universe, which is the basis of Newton's equation of motion,
 with the Earth coordinate system. This is used for satellite orbit calculation.\vskip.4em

$\bullet$\;  In seismology, it is necessary to quantify the two rotations \lq\lq{\it difference}\,\rq\rq\; to quantify how much the fault plane deviates from the reference plane.
\vskip1em 

For more purposes such as 
\vskip.5em

$\bullet$\; Computer graphics,
\vskip.4em

$\bullet$\; Aircraft design, spacecraft attitude control, dynamics such as aviation.
\vskip1em 

For more detailed comments, see e.g. [Agn], 2006. 

\vskip1.5em

{\bf A4.2.\; Time derivative of rotation in Rodriques expression}
\setcounter{equation}{1}

When we study rotations in a dynamical system, such as in aeronautical engineering, the rotation that depends on the time $t$ is treated. 
We give here an interesting compact formula for time derivative of rotation, which can be easily applied in many purposes.

 In case $\theta$ and ${\boldsymbol w}$ depend on the time $t$, put 
$Q(t):= R(\theta {\boldsymbol w})$, and let us study its time derivative \;$\dot{Q}(t)=\dfrac{dQ(t)}{dt}$. 

\vskip1em 

{\bf Theorem A4.1.}\; {\it 
 For \,${\boldsymbol w}\in {\boldsymbol B}_-\;(\subset{\boldsymbol H}_{\!-})$\, depending on $t$, we have \;$\dot{\boldsymbol{w}}\boldsymbol{w}=-\boldsymbol{w}\dot{\boldsymbol{w}}$\; and so \,$\dot{\boldsymbol{w}}\perp\boldsymbol{w}$\, in ${\boldsymbol H}_-$\,
 and \,$\dot{\boldsymbol{w}}\boldsymbol{w}=\dot{\boldsymbol{w}}\times\boldsymbol{w},\,\perp \dot{\boldsymbol{w}},\,\perp\boldsymbol{w}$ (mutually orthogonal). The time derivative $\dot{Q}(t)$ of rotation $Q(t)$ is given as 
\begin{eqnarray}
\label{2020-05-14-2}
&&
\dot{Q}(t)\,\boldsymbol{x} = Q(t)\Big(\big[\dot{\theta}\,\boldsymbol{w}+\sin \theta\;\dot{\boldsymbol{w}}+(1-\cos\theta)\,\dot{\boldsymbol{w}}\boldsymbol{w}\big]\times\boldsymbol{x} \Big),
\\
or 
\nonumber
&&
Q(t)^{-1}\dot{Q}(t)\,\boldsymbol{x} = \big[\dot{\theta}\,\boldsymbol{w}+\sin \theta\;\dot{\boldsymbol{w}}+(1-\cos\theta)\,\dot{\boldsymbol{w}}\boldsymbol{w}\big]\times\boldsymbol{x},
\qquad 
\end{eqnarray}
where $\boldsymbol{x}\in {\boldsymbol H}_-\;(\cong E^3)$ is a fixed vector, and for $\boldsymbol{a},\boldsymbol{b}\in {\boldsymbol H}_-$, \,$\boldsymbol{a}\times\boldsymbol{b}$ denotes the vector product in the 3-dimensional vector space ${\boldsymbol H}_-$.

}
\vskip1em
\setcounter{equation}{2}

{\it Proof.} \;Differentiate the both side of \;$\boldsymbol{w}^2 = \boldsymbol{w}\boldsymbol{w} = -1$\; with respect $t$, then  
\begin{eqnarray}
\label{2020-05-14-3}
\qquad
\dot{\boldsymbol{w}}\,\boldsymbol{w} +\boldsymbol{w}\,\dot{\boldsymbol{w}}\;=\;{\bf 0}\,\quad\therefore\; \boldsymbol{w} \perp \dot{\boldsymbol{w}}\;\;\mbox{\rm (in ${\boldsymbol H}_-$)}\quad\therefore\;  \dot{\boldsymbol{w}}\,\boldsymbol{w}=\dot{\boldsymbol{w}}\times\boldsymbol{w}.
\end{eqnarray}
\vskip-3em
\begin{gather}
\label{2020-05-14-4}
\nonumber 
\frac{d}{dt}\exp(\tfrac{1}{2}\theta\boldsymbol{w}) = \frac{d}{dt}\big(\cos(\tfrac{1}{2}\theta)+\sin(\tfrac{1}{2}\theta)\boldsymbol{w}\big)=
\\
\nonumber
\label{2020-05-14-5}
=\big(\!-\sin(\tfrac{1}{2}\theta)+\cos(\tfrac{1}{2}\theta)\boldsymbol{w}\big)\,\tfrac{1}{2}\dot{\theta}+ \sin(\tfrac{1}{2}\theta)\dot{\boldsymbol{w}}
=\exp(\tfrac{1}{2}\theta\boldsymbol{w})\,\tfrac{1}{2}\dot{\theta}\,\boldsymbol{w} + \sin(\tfrac{1}{2}\theta)\dot{\boldsymbol{w}}.
\quad 
\end{gather}
\setcounter{equation}{3}
Therefore 
\begin{eqnarray*} 
\quad
&&
\frac{d}{dt}Q(t)\boldsymbol{x}= \frac{d}{dt}\big\{\exp(\tfrac{1}{2}\theta\boldsymbol{w})\boldsymbol{x}\exp(-\tfrac{1}{2}\theta\boldsymbol{w})\big\}
\\
\nonumber 
&&
\quad
= \frac{d}{dt}\big(\exp(\tfrac{1}{2}\theta\boldsymbol{w})\big)\,\boldsymbol{x}\,\exp(-\tfrac{1}{2}\theta\boldsymbol{w})
+\exp(\tfrac{1}{2}\theta\boldsymbol{w})\,\boldsymbol{x}\,\frac{d}{dt}\big(\exp(-\tfrac{1}{2}\theta\boldsymbol{w})\big)
\\
\nonumber 
&&
\quad
= \exp(\tfrac{1}{2}\theta\boldsymbol{w})\big(\tfrac{1}{2}\dot{\theta}\big)\big(\boldsymbol{w}\boldsymbol{x}-\boldsymbol{x}\boldsymbol{w}\big)\exp(-\tfrac{1}{2}\theta\boldsymbol{w})+ 
\\
\nonumber
&&
\quad
+ \exp(\tfrac{1}{2}\theta\boldsymbol{w})\big\{\big(\cos(\tfrac{1}{2}\theta)-\sin(\tfrac{1}{2}\theta)\boldsymbol{w}\big)\, \sin(\tfrac{1}{2}\theta)\dot{\boldsymbol{w}}\boldsymbol{x} +
\\
\nonumber
&&
\quad
\hspace{14ex}
+\boldsymbol{x}\big(\!-\sin(\tfrac{1}{2}\theta)\dot{\boldsymbol{w}}\big)\big(\cos(\tfrac{1}{2}\theta)+\sin(\tfrac{1}{2}\theta)\boldsymbol{w}\big) \big\}\exp(-\tfrac{1}{2}\theta\boldsymbol{w})
\\
\nonumber 
&&
\quad
= \exp(\tfrac{1}{2}\theta\boldsymbol{w})\big\{\dot{\theta}\,\boldsymbol{w}\times\boldsymbol{x}+ \cos\big(\tfrac{1}{2}\theta\big)\sin\big(\tfrac{1}{2}\theta\big)(\dot{\boldsymbol{w}}\boldsymbol{x}-\boldsymbol{x}\dot{\boldsymbol{w}}\big)
\\
\nonumber
&&
\quad
\hspace{12ex}
- \sin^2\big(\tfrac{1}{2}\theta\big)\,\big(\boldsymbol{w}\,\dot{\boldsymbol{w}}\;\boldsymbol{x}\,+\,
\boldsymbol{x}\;\dot{\boldsymbol{w}}\;\boldsymbol{w}) \big\}\exp(-\tfrac{1}{2}\theta\boldsymbol{w})\qquad\mbox{\rm \big(use (\ref{2020-05-14-3})\big)}
\\
\nonumber 
&& 
\quad
= Q(t)\big\{\big[\dot{\theta}\,\boldsymbol{w}+\sin \theta\;\dot{\boldsymbol{w}}+(1-\cos\theta)\,\dot{\boldsymbol{w}}\boldsymbol{w}\big]\times\boldsymbol{x} \big\}.
\hspace{24ex}\Box
\end{eqnarray*}

\hspace*{24ex}{\bf \large (Omitted Below)}
\vskip2em

{\bf Acknowledgements.}\; The author would like to express heartful thanks to Prof.\;M.\,Duflo for documents of Rodrigues and to Prof.\;Kyo Nisiyama for documents of Hamilton.

{\small

\

}
\vskip2em

{\small  
Takeshi HIRAI, 22-8 Nakazaichi-Cho, Iwakura, Sakyo-Ku, Kyoto 606-0027, Japan; 

hirai.takeshi.24e@st.kyoto-u.ac.jp
}

\end{document}